



\def\updated{21 December 2004}

\count100= 46
\count101= 24
\message{: version 30nov04}

%





\font\titlefont=cmr17
\font\twelverm=cmr12
\font\ninerm=cmr9
\font\sevenrm=cmr7
\font\sixrm=cmr6
\font\fiverm=cmr5
\font\ninei=cmmi9
\font\seveni=cmmi7
\font\ninesy=cmsy9
\font\sevensy=cmsy9

\font\sixbf=cmbx6
\font\fivebf=cmbx5
\font\ninebf=cmbx9
\font\nineit=cmti9
\font\ninesl=cmsl9
\font\nineex=cmex9

\font\dfont=cmss10
\font\efont=cmti10

\def\ninepoint{\def\rm{\fam0\ninerm}
    \textfont0 = \ninerm
    \textfont1 = \ninei
    \textfont2 = \ninesy
    \textfont3 = \nineex
    \scriptfont0 = \sevenrm
    \scriptfont1 = \seveni
    \scriptfont2 = \sevensy
    \scriptscriptfont0 = \fiverm
    \scriptscriptfont1 = \fivei
    \scriptscriptfont2 = \fivesy
    \textfont\itfam=\nineit \def\it{\fam\itfam\nineit}
    \textfont\bffam=\ninebf \scriptfont\bffam=\sixbf
    \scriptscriptfont\bffam=\fivebf \def\bf{\fam\bffam\ninebf}
    \textfont\slfam=\ninesl \def\sl{\fam\slfam\ninesl}
    \baselineskip 10pt}


\def\CC{{\rm C\kern-.18cm\vrule width.6pt height 6pt depth-.2pt
\kern.18cm}}
\def\NN{{\mathop{{\rm I}\kern-.2em{\rm N}}\nolimits}}
\def\PP{{\mathop{{\rm I}\kern-.2em{\rm P}}\nolimits}}
\def\RR{{\mathop{{\rm I}\kern-.2em{\rm R}}\nolimits}}
\def\RRt{{\titlefont I}\kern-.2em{\titlefont R}}

\def\ZZ{{\mathop{{\rm Z}\kern-.28em{\rm Z}}\nolimits}}





\def\makebold#1{\mathord{\setbox0=\hbox{#1}%
       \copy0\kern-\wd0%
       \raise\dimen1\copy0\kern-\wd0%
       {\advance\dimen1 by \dimen1\raise\dimen1\copy0}\kern-\wd0%
       \kern\dimen0\raise\dimen1\copy0\kern-\wd0%
       {\advance\dimen1 by \dimen1\raise\dimen1\copy0}\kern-\wd0%
       \kern\dimen0\raise\dimen1\copy0\kern-\wd0%
       {\advance\dimen1 by \dimen1\raise\dimen1\copy0}\kern-\wd0%
       \kern\dimen0\raise\dimen1\copy0\kern-\wd0%
       \kern\dimen0\box0}}




\def \dword#1{{\dfont #1}}

\def \eword#1{{\efont #1}}

%


\def\norm#1{\Vert#1\Vert}

\def\dd{\,{\rm d}} 
\def\ii{{\rm i}}   

\def\frac#1#2{{#1 \over #2}}




\def\endproofsymbol{\makeblanksquare6{.4}}



\def\nopf{\medskip\goodbreak}

\def\makeblanksquare#1#2{
\dimen0=#1pt\advance\dimen0 by -#2pt
      \vrule height#1pt width#2pt depth0pt\kern-#2pt
      \vrule height#1pt width#1pt depth-\dimen0 \kern-#1pt
      \vrule height#2pt width#1pt depth0pt \kern-#2pt
      \vrule height#1pt width#2pt depth0pt
}

\magnification\magstep0

\hsize12.1truecm\vsize18.6truecm
\hoffset.8truein



\def\title#1{\toneormore#1||||:}
\def\titexp#1#2{\hbox{{\titlefont #1} \kern-.25em%
  \raise .90ex \hbox{\twelverm #2}}\/}
\def\titsub#1#2{\hbox{{\titlefont #1} \kern-.25em%
  \lower .60ex \hbox{\twelverm #2}}\/}

\def\author#1{\bigskip\bigskip\aoneormore#1||||:\smallskip\centerline{\updated}}

\def\abstract#1{\bigskip\bigskip\medskip%
 {\ninepoint
 \narrower{\bf Abstract.~}\rm#1\bigskip
 \printtochere}\starttoc\bigskip}

\def\toneormore#1||#2||#3:{\centerline{\titlefont #1}%
    \def\next{#2}\ifx\next\empty\else\medskip\toneormore#2||#3:\fi}
\def\aoneormore#1||#2||#3:{\centerline{\twelverm #1}%
    \def\next{#2}\ifx\next\empty\else\smallskip\aoneormore#2||#3:\fi}

\newwrite\toc\def\tocone{0}\def\tochalf{.5}\def\toctwo{1}
\def\printtochere{\immediate\closeout\toc{\inputifthere{\jobname.toc}}}
\def\starttoc{\immediate\openout\toc=\jobname.toc}
\def\nexttoc#1{{\let\folio=0\edef\next{\write\toc{#1}}\next}}

\def\tocline#1#2#3{\nexttoc{\line{\hskip\parindent\noexpand\noexpand\noexpand\rm\hskip#2truecm #1\hfill#3\hskip\parindent}}}

\def\footnoterule{\kern -3pt \hrule width 0truein \kern 2.6pt}
\def\leftheadline{\ifnum\pageno=\count100 \hfill%
  \else\rm\folio\hfil\it\shortauthor\fi}
\def\rightheadline{\ifnum\pageno=\count100 \hfill%
  \else\it\shorttitle\hfil\rm\folio\fi}

\nopagenumbers
\headline{\ifodd\pageno\rightheadline \else\leftheadline\fi}
\footline{\hfil}
\null\vskip 18pt
\centerline{}
\pageno=\count100
\count102=\count100
\advance\count102 by -1
\advance\count102 by \count101


\def\copyright{\hbox{{\twelverm o}\kern-.61em\raise .46ex\hbox{\fiverm c}}}

\insert\footins{\sixrm
\medskip
\baselineskip 8pt
\leftline{Surveys in Approximation Theory
  \hfill {\ninerm \the\pageno}}
\leftline{Volume 1, 2005.
pp.~\the\pageno--\the\count102.}
\leftline{Copyright \copyright\ 2005 Surveys in Approximation Theory.}
\leftline{ISSN x-x-x}
\leftline{All rights of reproduction in any form reserved.}
\smallskip
\par\allowbreak}


\def\sect#1{\startsect\edef\showsectno{\the\sectionno}\let\tocindent\tocone%
       \soneormore#1||||:\relax\medskip\noindent\ignorespaces}

\def\soneormore#1||#2||#3:{%
   \leftline{\bf\showsectno\hskip2em #1}
   \def\next{#2}%
   \ifx\next\empty\puttocline{\showsectno\ \ #1}{\folio}%
   \else\puttocline{\showsectno\ \ #1}{}\let\showsectno\skipsectno\let\tocindent\tochalf\soneormore#2||#3:\fi}

\def\puttocline#1#2{\tocline{#1}{\tocindent}{#2}}
\def\skipsectno{\setbox0=\hbox{\the\sectionno}\hskip\wd0}

\def\subsect#1{\formal{#1.}\let\tocindent\toctwo\puttocline{#1}{\folio}}
\def\formal#1{\bigskip{\bf #1}\hskip1em}

\def\fsubsect#1{\vskip-\baselineskip\subsect{#1}}

\newcount\sectionno\sectionno0
\def\presect{\the\sectionno.}
\newcount\subsectionno
\def\startsect{\ifx\empty\presect\else\restartnums\fi%
               \subsectionno0\global\advance\sectionno by 1\relax

               \goodbreak\bigskip\smallskip}
\def\startsubsect{\global\advance\subsectionno by 1\goodbreak\bigskip}


\def\figinbox#1(#2,#3)#4#5{\centerline{\vbox{\gridbox#2/#3/{
\ifshowfigname\point(0,0){#1}\fi 
  \point(0,0){\epsfxsize=#4truecm \epsfbox{\figsource #1}}#5}}}}

\input epsf
\def\figsource{}
\newif\ifshowfigname

\def\foneormore#1||#2||#3:{\centerline{\ninepoint\rm #1}%
    \def\next{#2}\ifx\next\empty\else\vskip0pt\foneormore#2||#3:\fi}


\def\gridbox#1/#2/#3{
\vbox to #1\gridunits{#3
\ifshowgrid\tickcount=0
  \loop\cgridw%
   \vbox to 0pt{\kern\tickcount \gridunits\hrule width#2\gridunits
       height\gridwidth\vss}
   \nointerlineskip \advance\tickcount by \tickskip
   \ifdim\tickcount pt<#1pt\repeat 
  \hbox to 0pt{\tickcount=0\hbox to 0pt{\tick#1/\hss}\advance\tickcount by \tickskip%
 \loop\ifdim\tickcount pt<#2pt\nexttick\tickskip#1/\advance\tickcount by \tickskip \repeat\hss}
\else \vbox to 0pt{\hrule width#2\gridunits height0pt\vss}
\fi\vfil}\vfil}

\def\ppoint#1#2(#3,#4)#5{\setbox0=\hbox{#5}
   \dimen0=\ht0\advance\dimen0 by\dp0\divide\dimen0 by-2
   \multiply\dimen0 by#1\advance\dimen0 by#3\gridunits
   \dimen1=\wd0\divide\dimen1 by-2\multiply\dimen1 by#2
   \advance\dimen1 by#4\gridunits\dpoint(\dimen0,\dimen1){#5}}

\def\gridunits{truecm}\newcount\tickskip\tickskip1\newcount\majortick\majortick5

\def\point(#1,#2)#3{\dpoint(#1\gridunits,#2\gridunits){#3}}
\def\dpoint(#1,#2)#3{\vbox to 0pt{\kern#1
   \hbox{\kern#2{#3}}\vss}\nointerlineskip}
\newcount\rmndr   
\def\rem#1#2{\rmndr=#1{}\divide\rmndr by#2{}%
\multiply\rmndr by-#2{}\advance\rmndr by #1}
\def\cgridw{\gridwidth\finegridw{}\rem\tickcount\majortick%
   \ifnum\rmndr=0{}\gridwidth\roughgridw\fi} 

\def\tick#1/{\cgridw\vrule width\gridwidth height0pt depth#1\gridunits}
\def\nexttick#1#2/{\hbox to#1\gridunits{\hfil\tick#2/}}
\newcount\tickcount
\newdimen\finegridw\finegridw0.4pt\newdimen\roughgridw\roughgridw1.6pt
\newdimen\gridwidth
\newif\ifshowgrid \showgridtrue





\def\label#1{%
  \ifsamelabel\global\samelabelfalse\else
  \ifmmode\global\advance\eqnum by 1
  \else\global\advance\labelnum by 1
  \fi\fi
  \edef\griff{label:#1}\edef\inhalt{\lastlabel}\definieres%
  \ifmmode\eqno(\inhalt)\else\inhalt\fi
  \ifdraft\ifmmode\rlap{\fiverm #1}\else\marginal{#1}\fi\fi}

\def\eqalignlabel#1{{\def\eqno{}\let\labelnum\eqnum\label{#1}}}

\def\labelplus#1#2{\def\labelsub{#2}\relax\label{#1}\def\labelsub{}}

\def\eqalignlabelplus#1#2{{\def\eqno{}\let\labelnum\eqnum\labelplus{#1}{#2}}}


\def\recall#1{\edef\griff{label:#1}\plazieres}

\newif\ifsamelabel
\def\labelsub{}
\def\lastlabel{\presect\ifmmode\the\eqnum\else\the\labelnum\fi\labelsub}
\def\nextlabel{{\ifmmode\advance\eqnum by 1\else\advance\labelnum by 1\fi\lastlabel}}


\newcount\blackmarks\blackmarks0
\newcount\eqnum
\newcount\labelnum
\def\restartnums{\eqnum0\labelnum0}
\def\singlecount{\let\labelnum\eqnum}

 \newread\testfl
 \def\inputifthere#1{\immediate\openin\testfl=#1
    \ifeof\testfl\message{(#1 does not yet exist)}
    \else\input#1\fi\closein\testfl}

 \inputifthere{\jobname.aux}
 \newwrite\aux
 \immediate\openout\aux=\jobname.aux

\def\plazieres{\expandafter\ifx\csname\griff\endcsname\relax%
  \xdef\esfehlt{\griff}\blackmark\else{\csname\griff\endcsname}\fi}

\def\definieres{\expandafter\xdef\csname\griff\endcsname{\inhalt}%
 \def\blankkk{ }\expandafter\immediate\write\aux{%
 \string\expandafter\def\string\csname%
 \blankkk\griff\string\endcsname{\inhalt}}}

\def\blackmark{\ifnum\blackmarks=0\global\blackmarks=1%
 \write16{============================================================}%
 \write16{Some forward reference is not yet defined. Re-TeX this file!}%
 \write16{============================================================}%
 \fi\immediate\write16{undefined forward reference: \esfehlt}%
 {\vrule height10pt width2pt depth2pt}\esfehlt%
 {\vrule height10pt width2pt depth2pt}}

\def\marginal#1{\strut\vadjust{\kern-\strutdepth%
\vtop to \strutdepth{\baselineskip\strutdepth\vss\llap{\fiverm#1\ }\null}}}
\def\strutdepth{\dp\strutbox}


\newif\ifdraft

\newcount\hour\newcount\minutes
\def\draft{\drafttrue
\def\comment##1{{\bf comment: ##1}}
\headline={\sevenrm \hfill\ifx\filenamed\undefined\jobname\else\filenamed\fi%
(.tex) (as of \ifx\updated\undefined???\else\updated\fi)
 \TeX'ed at {\hour\time\divide\hour by 60{}%
\minutes\hour\multiply\minutes by 60{}%
\advance\time by -\minutes
\the\hour:\ifnum\time<10{}0\fi\the\time\  on \today\hfill}}
}

\def\today{\number\day\space%
\ifcase\month\or January\or February\or March\or April\or May\or June\or
 July\or August\or September\or October\or November\or December\fi%
\space\number\year}



\def\References{\goodbreak\bigskip\centerline{\bf References}%
   \tocline{\skipsectno\ \  References}{\tocone}{\folio}%
   \bigskip\frenchspacing}

\def\bibitem{\smallskip\noindent}


\gdef\formfirstauthor{\the\lastname, \the\firstname}
\gdef\formnextauthor{, \the\firstname\the\lastname}
\gdef\formotherauthor{ and \the\firstname\the\lastname}
\gdef\formlastauthor{,\formotherauthor}

\gdef\formB{\the\au\ [\yr] ``\the\ti'', \the\pb, \the\pl. \setcitelabel}
\gdef\formD{\the\au\ [\yr] ``\the\ti'', dissertation, \the\pl. \setcitelabel}
\gdef\formJ{\the\au\ [\yr] \the\ti, {\sl\the\jr}\ifx\vl\empty%
\else\ {\bf\vl}\fi, \pp. \setcitelabel}
\gdef\formP{\the\au\ [\yr] \the\ti, in {\sl\the\tit},
\getfirstchar\aut\ifx\firstchar\unknownx\else\the\aut, ed\edsop, \fi
\getfirstchar\pub\ifx\firstchar\unknownx\else\the\pub, \fi \the\pl, \pp. \setcitelabel}
\gdef\formR{\the\au\ [\yr] \the\ti\ifx\is\empty\else, \is\fi. \setcitelabel}

\newtoks\lastname
\newtoks\firstname
\newtoks\au
\newtoks\aut
\newtoks\ti
\newtoks\tit
\newtoks\pb
\newtoks\pub
\newtoks\pl
\newtoks\jr

\newtoks\rhlau

\def\setcitelabel{\edef\griff{cit\rh}\edef\inhalt{\the\rhlau\ \yr}\definieres}
\def\setcitelabel{}

\def\getfirstchar#1{\edef\theword{\the#1}\expandafter\getit\theword:}
\def\getit#1#2:{\def\firstchar{#1}}
\def\unknownx{x}

\newif\ifonesofar
\def\concat#1{\edef\audef{{#1}}\au=\audef}
\def\decodeauthor#1, #2,#3;{\lastname={#1}\firstname={#2}%
\concat{\formfirstauthor}\onesofartrue%
\def\morerhlau{}%
\def\next{#3}\ifx\next\empty\else\def\morerhlau{ et al.}\decodemoreauthor#3;\fi
\edef\morerhlauu{{\the\lastname\morerhlau}}\rhlau=\morerhlauu}
\def\decodemoreauthor#1, #2,#3;{\lastname={#1}\firstname={#2}%
\def\next{#3}\ifx\next\empty\let\formaut=\formlastauthor%
\ifonesofar\ifx\formotherauthor\undefined\else\let\formaut=\formotherauthor%
\fi\fi\concat{\the\au\formaut}%
\else\onesofarfalse\concat{\the\au\formnextauthor}\decodemoreauthor#3;\fi}

\def\refproc #1(#2; #3; {\decodeproc#2; \xdef\yr{#3}}
\def\decodeproc#1), #2 (ed#3.),#4 (#5); {%
 \global\tit={#1}\global\aut={#2}\xdef\edsop{#3}\global
 \pub={#4}\global\pl={#5}}


\def\refB #1; #2; #3 (#4); #5; {\decodeauthor#1,;%
   \ti={#2}\pb={#3}\pl={#4}\def\yr{#5}\bibitem\formB}

\def\refD #1; #2; #3; #4; {\decodeauthor#1,;%
   \ti={#2}\pl={#3}\def\yr{#4}\bibitem\formD}

\def\refJ #1; #2; #3; #4; #5; #6; {\decodeauthor#1,;%
    \ti={#2}\jr={#3}\def\vl{#4}\def\yr{#5}\def\pp{#6}\bibitem\formJ}

\def\refP #1; #2; #3; #4; {{\global\aut={\vrule height15pt width15pt depth0pt}%
 \global\tit={{\bf the specified proceedings does not exist in our files}}%
 \xdef\edsop{}\global\pub={}\def#3{}\input proceed }\decodeauthor#1,;%
        \ti={#2}\def\pp{#4}\bibitem\formP}

\def\refQ #1; #2; (#3; #4; #5; {\decodeproc#3; \decodeauthor#1,;%
   \ti={#2}\def\yr{#4}\def\pp{#5}\bibitem\formP}

\def\refR #1; #2; #3; #4; {\decodeauthor#1,;%
         \ti={#2}\def\is{#3}\def\yr{#4}\bibitem\formR}

\singlecount   
\def\presect{} 


\title{Divided Differences}
\author{Carl de Boor}

\def\shorttitle{Divided Differences}
\def\shortauthor{C. de Boor}


\def\divdif{\mathord{\kern.43em\vrule width.6pt height5.6pt depth-.28pt\kern-.42em\Delta}}
\def\dvd#1{\divdif(#1)}
\newcount\excount\excount0
\def\nextex{\global\advance\excount by 1{}\formal{Example \the\excount.}}
\def\ootpii{{1\over2\pi\ii}}
\def\nwt#1{w_{#1}}
\def\Nwt#1{W_{#1}}

\def\hatc{\widehat{c}}
\def\hatt{\widehat{t}}
\def\Matrix#1{\left[\matrix{#1}\right]}
\def\braket#1{\hbox{$[\![$}#1\hbox{$]\!]$}}
\let\bs\backslash
\def\dist{\mathop{\rm dist}{}}
\def\ran{\mathop{\rm ran}\nolimits}
\def\Null{\mathop{\rm null}\nolimits}
\def\fromto{\mathbin{\ldotp\ldotp}} 
\let\gD\Delta
\let\gd\delta
\let\ga\alpha
\let\gz\zeta
\let\gl\lambda
\let\gs\sigma
\let\gt\tau
\def\Fn{\FF^n}
\def\FF{{{\rm I}\kern-.16em {\rm F}}}
\def\inv#1{#1^{-1}}
\def\id#1{{\rm id}{}_{#1}}
\def\nextline{\hfill\break\noindent}
\let\newline\nextline
\def\newpar{\par}
\def\Iff{\hskip1em\Longleftrightarrow\hskip1em}

\let\proof\pf
\def\endproof{\hbox{}~\hfill\endproofsymbol\nopf}
\let\endexample\endproof

\def\rhl#1{\rahel#1::/}
\def\rahel#1:#2:#3/{\edef\rh{#1}\def\next{#2}\ifx\next\empty\edef\rl{#1}%
\else\edef\rl{#2}\fi%
\edef\griff{cit\rh}\edef\inhalt{\rl}\definieres}

\def\cite#1{\excite#1::/}
\def\excite#1:#2:#3/{\edef\griff{cit#1}%
\formcite{\plazieres\def\next{#2}\ifx\next\empty\else\formciteaddl{#2}\fi}}
\def\formcite#1{{\bf[}#1{\bf]}}
\def\formciteaddl#1{: #1}
\bigskip
{\ninepoint
\rightline{\it``Est enim fere ex pulcherrimis qu\ae\ solvere
desiderem.''}
\rightline{\rm
(It is among the most beautiful I could desire to solve.) \cite{N76}}}
\vskip-\baselineskip

\abstract{Starting with a novel definition of divided differences, this essay
derives and discusses the basic properties of, and facts about,
(univariate) divided differences.}


\sect{Introduction and basic facts}
While there are several ways to think of divided differences, including the one
suggested by their
very name, the most efficient way is as the coefficients in a Newton form.
This form provides an efficient representation of Hermite interpolants.

Let $\Pi\subset(\FF\to\FF)$ be the linear space of polynomials in one real 
($\FF=\RR$) or complex ($\FF=\CC$) variable, and let $\Pi_{<n}$ denote the 
subspace of all polynomials of degree $<n$.
The \dword{Newton form} of $p\in\Pi$ with respect to the sequence $t=(t_1,t_2,\ldots)$ of
\dword{centers} $t_j$ is its expansion
$$
p =: \sum_{j=1}^\infty \nwt{j-1,t}\,c(j)
\label{newtonform}
$$
in terms of the \dword{Newton polynomial}s
$$
\nwt i\;:=\;\nwt{i,t}\;:=\;(\cdot-t_1)\cdots(\cdot-t_i),\quad i=0,1,\ldots\;.
\label{defnwt}
$$
Each $p\in\Pi$ does, indeed, have exactly one such expansion for any
given $t$ since $\deg\nwt{j,t}=j$, all $j$, hence $(w_{j-1,t}: j\in\NN)$ is a
\dword{graded basis for $\Pi$} in the sense that, for each $n$, $(w_{j-1,t}:
j=1{:}n)$ is a basis for $\Pi_{<n}$.

In other words, the column map
$$\Nwt t:\FF_0^\NN\to\Pi: c \mapsto \sum_{j=1}^\infty \nwt{j-1,t}c(j)
\label{defVt}$$
(from the space $\FF_0^\NN$ of scalar sequences with finitely many nonzero
entries to the space $\Pi$)
is 1-1 and onto, hence invertible.
In particular, for each $n\in\NN$, the coefficient $c(n)$ in the Newton form
(\recall{newtonform}) for $p$ depends linearly on $p$, i.e.,
{\sl $p\mapsto c(n) = (\inv{\Nwt t}p)(n)$ is a
well-defined linear functional on $\Pi$, and vanishes on $\Pi_{<n-1}$.}
More than that,
since all the (finitely many nontrivial) terms in (\recall{newtonform}) with $j>n$ have
$w_{n,t}$ as a factor, we can write
$$
p\;=\; p_n + \nwt{n,t} q_n,
\label{defpnqn}$$
with $q_n$ a polynomial we will look at later (in Example 6), and with
$$p_n\;:=\;\sum_{j=1}^n\nwt {j-1,t}c(j)$$
a polynomial of degree $<n$. This makes $p_n$
necessarily the remainder left by the division of $p$ by $\nwt{n,t}$, hence
well-defined for every $n$, hence, by induction, we obtain another proof that
the Newton form (\recall{newtonform}) itself is well-defined.

In particular, $p_n$ depends only on $p$ and on
$$t_{1:n}:= (t_1,\ldots,t_n),$$
therefore the same is true of its leading coefficient, $c(n)$.
This is reflected in the (implicit) definition
$$
p\;=:\; \sum_{j=1}^\infty \nwt{j-1,t} \dvd{t_{1:j}}p,\quad p\in\Pi,
\label{idefdvd}
$$
in which the coefficient $c(j)$ in the Newton form (\recall{newtonform}) for $p$ is
denoted
$$\dvd{t_{1:j}}p\;=\;\dvd{t_1,\ldots,t_j}p\;:=\;(\inv{(\Nwt t)}p)(j)
\label{defdvd}
$$
and called the \dword{divided difference of
$p$ at $t_1,\ldots,t_j$}. It is also called a divided difference of
\dword{order} $j-1$, and the reason for all this terminology will
be made clear in a moment.

Since $\Nwt t$ is a continuous function of $t$, so is $\inv{\Nwt t}$, hence so
is $\dvd{t_{1:j}}$ (see Proposition \recall{propcont} for proof details).
Further, since $w_{j,t}$ is symmetric in $t_1,\ldots,t_j$,
so is $\dvd{t_{1:j}}$. Also, $\dvd{t_{1:j}}\perp\Pi_{<j}$ (as mentioned before).

In more practical terms, we have

\proclaim Proposition \label{prophermiteconds}. The sum
$$p_n\;=\;\sum_{j=1}^n\nwt{j-1,t}\dvd{t_{1:j}}p$$
of the first $n$ terms in the Newton form (\recall{newtonform}) for $p$ is the
\dword{Hermite interpolant to $p$ at $t_{1:n}$},
i.e., the unique
polynomial $r$ of degree $<n$ that \dword{agrees with $p$ at $t_{1:n}$}
in the sense that
$$D^i r(z) = D^i p(z),\quad 0\le i<\mu_z:=\#\{j\in[1\fromto n]: t_j=z\},\quad
z\in \FF.
\label{hermiteconds}
$$

\proof
One readily verifies by induction on the nonnegative integer $\mu$ that, for
any $z\in\FF$, any polynomial $f$ \dword{vanishes $\mu$-fold at $z$} iff
$f$ has $(\cdot-z)^\mu$ as a factor, i.e.,
$$D^if(z) = 0 \hbox{\ \ for\ \ } i=0,\ldots,\mu-1 \Iff f \in
(\cdot-z)^\mu\Pi.\label{mufoldfactor}$$

Since $p-p_n = \nwt{n,t}q_n$, this implies that $r=p_n$ does,
indeed satisfy (\recall{hermiteconds}).

Also, $p_n$ is the only such polynomial since, by (\recall{mufoldfactor}), for any
polynomial $r$ satisfying (\recall{hermiteconds}), the difference $p_n-r$ must have
$\nwt n$ as a factor and, if $r$ is of degree $<n$, then this is possible only
when $r=p_n$.
\endproof

\nextex  For $n=1$, we get that
$$\dvd{t_1}:p\mapsto p(t_1),$$
i.e., $\dvd{\gt}$ can serve as a (nonstandard) notation for the linear
functional of evaluation at $\gt$.\endexample

\nextex  For $n=2$, $p_n$ is the polynomial of degree $<2$ that
matches $p$ at $t_{1:2}$. If $t_1\not=t_2$, then we know $p_2$ to be
writable in `point-slope form' as
$$p_2 = p(t_1) + (\cdot-t_1){p(t_2)-p(t_1)\over t_2-t_1},$$
while if $t_1=t_2$, then we know $p_2$ to be
$$p_2 = p(t_1) + (\cdot-t_1)Dp(t_1).$$
Hence, altogether,
$$
\dvd{t_{1:2}}p\;=\;
\cases{ {\displaystyle p(t_2)-p(t_1)\over \displaystyle\mathstrut t_2-t_1},&$t_1\not=t_2$;\cr
\noalign{\vskip4pt}
         Dp(t_1),&otherwise.\cr}\label{firstdvd}
     $$
Thus, for $t_1\not=t_2$, $\dvd{t_{1:2}}$ is
a quotient of differences, i.e., a \eword{divided difference}.\endexample

\nextex   Directly from the definition of the divided difference,
$$\dvd{t_{1:j}}\nwt{i-1,t} = \gd_{ji},
\label{biorthogonal}
$$
therefore (remembering that $\dvd{t_{1:j}}\perp\Pi_{<j-1}$)
$$
\dvd{t_{1:j}}()^{j-1}\;=\;1,
\label{dvdofpower}
$$
with
$$()^k:\FF\to\FF: z\mapsto z^k$$
a handy if nonstandard notation for the power functions.
\endexample

\nextex  If $t$ is a constant sequence, $t=(\gt,\gt,\ldots)$ say,
then $$\nwt{j,(\gt,\gt,\ldots)} = (\cdot-\gt)^{j},$$
hence the \dword{Taylor expansion}
$$
p = \sum_{n=0}^\infty (\cdot-\gt)^n D^np(\gt)/n!\label{taylorexp}
$$
is the Newton form for the polynomial $p$ with respect to the sequence
$(\gt,\gt,\ldots)$.
Therefore,
$$\dvd{\gt^{[n+1]}}p := \dvd{\underbrace{\gt,\ldots,\gt}_{n+1\,\rm terms}}p =
D^np(\gt)/n!,\quad n=0,1,\ldots\,.
\label{dvdallthesame}
$$
\endexample

\nextex  If $\ell: t\mapsto at+b$, then
$(\ell(z)-\ell(t_i)) = a(z - t_i)$, hence
$$a^{n-1}\dvd{\ell(t_{1:n})}p = \dvd{t_{1{:}n}}(p\circ\ell).
\label{changevar}
$$
\endexample

\nextex  Consider the polynomial $q_n$ introduced in (\recall{defpnqn}):
$$
p = p_n + \nwt{n,t}q_n.$$
Since
$p(t_{n+1}) = p_{n+1}(t_{n+1})$
and
$p_{n+1} = p_n + w_{n,t}\dvd{t_{1:n+1}}p$,
we have
$$
\nwt{n,t}(t_{n+1})q_n(t_{n+1}) = \nwt{n,t}(t_{n+1})\dvd{t_{1:n+1}}p,$$
therefore
$$q_n(t_{n+1}) = \dvd{t_{1:n+1}}p,$$
at least for any $t_{n+1}$ for which $\nwt{n,t}(t_{n+1})\not=0$,
hence for every
$t_{n+1}\in\FF$, by the continuity of $q_n$, and the continuity of
$\dvd{t_{1{:}n},\cdot}p$, i.e., of
$\dvd{t_{1{:}n+1}}p$ as a function of $t_{n+1}$.
It follows that
$$q_n=\dvd{t_{1:n},\cdot}p$$
and
$$p = p_n + \nwt{n,t}\dvd{t_{1:n},\cdot}p,
\label{withremainder}
$$
the \dword{standard error formula for Hermite interpolation}.
More than that, by the very definition, (\recall{defpnqn}), of $q_n$, we now know that
$$\dvd{t_{1:n},\cdot}p
\;=\;q_n = (p-p_n)/\nwt{n,t} =
\sum_{j>n}{\nwt{j-1,t}\over\nwt{n,t}}\dvd{t_{1:j}}p,
\label{dvdasfunction}$$
and we recognize the sum here as a Newton form with respect to the
sequence $(t_j: j>n)$. This provides us with the following
\dword{basic divided difference identity}:
$$
\dvd{t_{n+1:j}}\dvd{t_{1:n},\cdot}
= \dvd{t_{1:j}},\quad j>n.
\label{compofdvds}
$$
\endexample

For the special case $n=j-2$, the basic divided difference identity,
(\recall{compofdvds}), reads
$$
\dvd{t_{j-1:j}}\dvd{t_{1:j-2},\cdot}=\dvd{t_{1:j}},$$
or, perhaps more suggestively,
$$
\dvd{t_{j-1},\cdot}\dvd{t_{1:j-2},\cdot}=\dvd{t_{1:j-1},\cdot},$$
hence, by induction,
$$
\dvd{t_{j-1},\cdot}\dvd{t_{j-2},\cdot}\,\cdots\dvd{t_1,\cdot}
\;=\; \dvd{t_{1:j-1},\cdot}.
\label{dvdisdvd}$$
In other words, $\dvd{t_{1:j}}$ is obtainable by forming difference quotients
$j-1$ times. This explains our calling $\dvd{t_{1:j}}$ a
`divided difference of order $j-1$'.

%
%
%
%
%
%

\sect{Continuity and smoothness}
{\sl The column map
$$\Nwt t:\FF_0^\NN\to\Pi: c\mapsto \sum_{j=1}^\infty \nwt{j-1,t}c(j)
$$
{\rm introduced in (\recall{defVt})} is continuous as a function of $t$, hence
so is its inverse}, as follows directly from
the identity
$$\inv{A}-\inv{B}= \inv{A}(B-A)\inv{B},\label{diffinv}$$
valid for any two invertible maps $A$, $B$ (with the same domain and target).
Therefore, also each $\dvd{t_{1:j}}$ is a continuous function of $t$,
all of this in the pointwise sense. Here is the formal statement and its proof.
\proclaim Proposition \label{propcont}. For any $p\in\Pi$,
$$\lim_{s\to t} (\dvd{s_{1:j}}p: j\in\NN) = (\dvd{t_{1:j}}p: j\in\NN).
$$

\proof
Let $p\in\Pi_{<n}$. Then $t\mapsto\dvd{t_{1{:}k}}p = 0$ for
$k>n$, hence trivially continuous. As for $k\le n$, let
$$\Nwt{t,n}:=\Fn\to\Pi_{<n}: c\mapsto \sum_{j=1}^n \nwt{j-1,t}c(j)$$
be the restriction of $\Nwt t$ to $\Fn$, as a linear map to $\Pi_{<n}$.
Then, in whatever norms we might choose on $\Fn$ and $\Pi_{<n}$, $\Nwt{t,n}$ is
bounded and invertible, hence boundedly invertible uniformly in
$t_{1:n}$ as long as $t_{1:n}$ lies in some
bounded set.
Therefore, with (\recall{diffinv}), since $\lim_{s\to t}\Nwt{s,n} = \Nwt{t,n}$, also
$$\lim_{s\to t} (\dvd{s_{1:j}}p: j=1{:}n) =
\inv{(\Nwt{t,n})}p =
(\dvd{t_{1:j}}p: j=1{:}n).
$$
\endproof

This continuity is very useful. For example, it implies that it is usually
sufficient to check a proposed
divided difference identity by checking it only for pairwise distinct arguments.

As another example, we used  the continuity earlier (in Example 6) to prove
that $\dvd{t_{1{:}n},\cdot}p$
is a polynomial. This implies that $\dvd{t_{1{:}n},\cdot}p$ is differentiable,
and, with that, (\recall{compofdvds}) and (\recall{dvdallthesame})  even provide the
following formula for
the derivatives.
\proclaim Proposition \label{propdiff}.
$$D^k\dvd{t_{1:j},\cdot}p = k!\dvd{t_{1:j},[\cdot]^{k+1}}p,\quad k\in \NN.$$

\sect{Refinement}
Already Cauchy \cite{Ca} had occasion to use the simplest nontrivial case of
the following fact.

\proclaim Proposition \label{proprefine}.
For any $n$-sequence $t$ and any $1\le \gs(1)<\cdots<\gs(k)\le n$,
$$
\dvd{t_{\gs(1:k)}} =
\sum_{j=\gs(1)-1}^{\gs(k)-k}\ga(j)\dvd{t_{j+1:j+k}},$$
with $\ga = \ga_{t,\gs}$ positive in case $t$ is strictly increasing.

\proof Since $\dvd{t_{1:n}}$ is symmetric in the $t_j$, (\recall{compofdvds})
implies
$$
(t_n-t_1)(\dvd{t_{1:n\bs m}}-\dvd{t_{2:n}})
\;=\;
(t_1-t_m)(\dvd{t_{2:n}}-\dvd{t_{1:n-1}}),
$$
with
$$t_{1{:}n\bs m}:=
t_{1{:}m-1,m+1{:}n}:=(t_1,\ldots,t_{m-1},t_{m+1},\ldots,t_n).$$
On rearranging the terms, we get
$$(t_n-t_1)\dvd{t_{1:n\bs m}}
\;=\;
(t_n-t_m)\dvd{t_{2:n}}+(t_m-t_1)\dvd{t_{1:n-1}},$$
and this proves the assertion for the special case $k=n-1$, and even gives an
explicit formula for $\ga$ in this case.

From this, the general case follows by induction on $n-k$, with $\ga$
computable as a convolution of sequences which, by induction, are positive in
case $t$ is strictly increasing (since this is then trivially so for
$k=n-1$), hence then $\ga$ itself is positive.  \endproof

My earliest reference for the general case is \cite{P33}.

\sect{Divided difference of a product: Leibniz, Opitz}
\def\Piln{\Pi_{<n}}
The map
$$P := P_{n,t}:\Pi \to \Pi: p\mapsto p_n,$$
of Hermite interpolation at $t_{1:n}$,
is the linear projector $P$ on $\Pi$ with
$$\ran P = \Piln,\quad \ran(\id{}-P) = \Null P = w_{t,n}\Pi.$$
In particular, the nullspace of $P$ is an \eword{ideal} if, as we may, we think
of $\Pi$ as a ring, namely the ring with multiplication defined pointwise,
$$(pg)(z) := p(z)g(z),\quad z\in\FF.$$
In other words, the nullspace of $P$ is a linear
subspace closed also
under pointwise multiplication. This latter fact is (see \cite{B03b})
equivalent to the identity
$$P(pq) = P(p(Pq)),\quad p,q\in\Pi.\label{charidealproj}$$

For $p\in\Pi$, consider the map
$$M_p:\Piln\to\Piln:f\mapsto P(pf).$$
Then $M_p$ is evidently linear and, also evidently, so is the resulting map
$$M:\Pi\to L(\Piln): p\mapsto M_p$$
on $\Pi$ to the space of linear maps on $\Piln$. More than that, since, by
(\recall{charidealproj}),
$$M_{pq}f = P(pqf) = P(pP(qf)) = M_pM_qf,\quad p\in\Piln,\; p,q\in\Pi,$$
$M$ is a ring homomorphism, from the ring $\Pi$
into the ring $L(\Piln)$ in which composition serves as multiplication.
The latter ring is
well known not to be commutative while, evidently, $\ran M$ is a commutative
subring.

It follows, in particular, that
$$M_p = p(M_{()^1}),\quad p\in \Pi,$$
hence
$$\widehat{M_p} = p(\widehat{M_{()^1}}),\quad p\in \Pi,$$
for the \eword{matrix representation}
$$\widehat{M_p} := VM_p\inv{V}$$
of $M_p$ with respect to any particular basis $V$ of $\Piln$.
Look, in particular, at the matrix representation with respect to the Newton
basis
$$
V := [w_{j-1,t}: j=1{:}n]
$$
for $\Piln$.
Since
$$ ()^1w_{j-1,t} = t_jw_{j-1,t} + w_{j,t},\quad j=1,2,\ldots,$$
therefore evidently
$$M_{()^1}w_{j-1,t} = P(()^1w_{j-1,t})= t_jw_{j-1,t} + (1-\gd_{j,n}) w_{j,t},
\quad j=1{:}n.$$
Consequently, the matrix representation for $M_{()^1}$ with respect to the
Newton basis $V$ is the bidiagonal matrix
$$
\widehat{M_{()^1}}\;=\;A_{n,t}\;:=\;
\Matrix{t_1&   &   &   & \cr
         1 &t_2&   &   & \cr
           & 1 &t_3&&\cr
                  &   & \ddots&\ddots&\cr
                  &   &   &  1   & t_n \cr}.
$$
On the other hand, for any $p\in\Pi$ and $j=1{:}n$,
$$\left(\sum_{i=j}^n (w_{i-1,t}/w_{j-1,t})\dvd{t_{j:i}}p\right)w_{j-1,t}$$
is a polynomial of degree $<n$ and, for pairwise distinct $t_i$,
it agrees with $p w_{j-1,t}$ at $t_{1:n}$ since the sum describes the polynomial
of degree $\le n-j$ that matches $p$ at $t_{j:n}$ while both functions vanish at
$t_{1:j-1}$. Consequently, with the convenient agreement that
$$\dvd{t_{j:i}}:=0, \qquad j>i,$$
we conclude that
$$P(pw_{j-1,t}) = \sum_{i=1}^n w_{i-1,t}\dvd{t_{j:i}}p,\quad j=1{:}n,$$
at least when the $t_i$ are pairwise distinct.
In other words, the $j$th column of the matrix $\widehat{M_p} = \inv{V}M_pV$
(which represents $M_p$ with respect to the Newton basis $V$ for $\Piln$)
has the entries
$$(\dvd{t_{j:i}}p: i=1{:}n) =
(0,\ldots,0,p(t_j),\dvd{t_j,t_{j+1}}p,\ldots,\dvd{t_{j{:}n}}p).$$
By the continuity of the divided difference (see Proposition \recall{propcont}),
this implies
\proclaim Proposition \label{propopitz}: Opitz formula. For any $p\in\Pi$,
$$p(A_{n,t})\;=\;(\dvd{t_{j:i}}p: i,j=1{:}n).\label{opitzformula}
$$

The remarkable identity (\recall{opitzformula}) is due to G. Opitz; see \cite{O} which
records a talk announced but not delivered. Opitz calls the matrices
$p(A_{n,t})$ \eword{Steigungsmatrizen} (`difference-quotient matrices').
Surprisingly, Opitz explicitly
excludes the possibility that some of the $t_j$ might coincide.
\cite{BR} ascribe (\recall{opitzformula}) to Sylvester, but I have been unable to
locate anything like this formula in Sylvester's collected works.

\nextex  For the monomial $()^k$, Opitz' formula gives
$$\dvd{t_{1:n}}()^k = (A_{n,t})^k(n,1)
= \sum_{\nu\in\{1{:}n\}^k} A_{n,t}(n,\nu_k) A_{n,t}(\nu_k,\nu_{k-1})\cdots A_{n,t}(\nu_1,1),$$
and, since $A_{n,t}$ is bidiagonal, the $\nu$th summand is zero unless the
sequence
$(1,\nu_1,\ldots,\nu_k,n)$
is increasing, with any strict increase no bigger than 1, in which case the
summand equals $t^\ga$, with $\ga_j-1$ the multiplicity with which $j$ appears
in the sequence $\nu$, $j=1{:}n$. This confirms that
$\dvd{t_{1:n}}()^k = 0$ for $k<n-1$ and proves that
$$
\dvd{t_{1:n}}()^k = \sum_{|\ga|=k-n-1}t^\ga,\quad k\ge n-1.
\label{dvdofpower}
$$
\endexample

My first reference for (\recall{dvdofpower}) is \cite{St27:p.19f}.
To be sure, once (\recall{dvdofpower}) is known, it is easily verified by
induction, using the Leibniz formula, to be derived next.

Since, for any square matrix $A$ and any polynomials $p$ and $q$,
$$(pq)(A) = p(A)q(A),$$
it follows, in particular, that
$$\dvd{t_{1:n}}(pq) = \widehat{M_{pq}}(n,1) =
\widehat{M_p}(n,:)\widehat{M_q}(:,1),$$
hence
\proclaim Corollary \label{corleibniz}: Leibniz formula. For any $p,q\in\Pi$,
$$\dvd{t_{1:n}}(pq) = \sum_{j=1:n}\dvd{t_{j:n}}p\;\dvd{t_{1:j}}q.
\label{leibnizformula}$$

On the other hand, the Leibniz formula implies that, for any $p,q\in\Pi$,
$$(\dvd{t_{j:i}}p: i,j=1{:}n)
  (\dvd{t_{j:i}}q: i,j=1{:}n)
  \;=\;
  (\dvd{t_{j:i}}(pq): i,j=1{:}n),$$
hence, that, for any $p\in\Pi$,
$$p((\dvd{t_{j:i}}()^1: i,j=1{:}n))
\;=\;
(\dvd{t_{j:i}}p: i,j=1{:}n).$$
In other words, we can also view Opitz' formula as a corollary to Leibniz'
formula.

My first reference for the Leibniz formula is
\cite{P33}, though Steffensen later devotes
an entire paper, \cite{St39}, to it and this has become the standard reference
for it despite the fact that Popoviciu, in response, wrote
his own overview of divided differences, \cite{P40}, trying, in vain, to
correct the record.

The (obvious) name `Leibniz formula' for it appears first in \cite{B72}.

Induction on $m$ proves the following

\proclaim Corollary \label{genleibniz}: General Leibniz formula.
For $f:\FF^m\to\FF$,
$$
\dvd{t_1,\ldots,t_k}f(\cdot,\ldots,\cdot) =
\sum_{1=i(1)\le\cdots\le i(m)=k}\left(\otimes_{j=1}^m
\dvd{t_{i(j-1)},\ldots,t_{i(j)}}\right)f.
$$


\sect{Construction of Newton form via a divided difference table}

\proclaim Divided difference table. Assume that the sequence $(t_1,\ldots,t_n)$
has all its multiplicities (if any) \dword{clustered}, meaning
that, for any $i<j$, $t_i = t_j$ implies that $t_i=t_{i+1}=\cdots=t_j$.
Then, by (\recall{compofdvds}) and (\recall{dvdallthesame}),
$$\dvd{t_{i:j}}p\;=\;\cases{
{{\displaystyle\mathstrut\dvd{t_{i+1:j}}p-\dvd{t_{i:j-1}}p}\over
{\displaystyle\mathstrut t_j-t_i}},&$t_i\not=t_j$;\cr
\noalign{\vskip3pt}
D^{j-i}p(t_i)/(j-i)! &otherwise,\cr}\quad 1\le i\le j\le n.
$$
Hence, it is possible to fill in all the entries in the
\dword{divided difference table}
$$
\matrix{
\dvd{t_1}p&&&&\cr
&\dvd{t_{1:2}}p&&&\cr
\dvd{t_2}p&&\dvd{t_{1:3}}p&&\cr
&\dvd{t_{2:3}}p&&\ \ \ \cdot\ \ \ &\cr
\dvd{t_3}p&&\cdot&&\dvd{t_{1:n-1}}p\cr
&\cdot&&\ \ \ \cdot\ \ \ &&\dvd{t_{1:n}}p\cr
\cdot&&\dvd{t_{n-3:n-1}}p&&\dvd{t_{2:n}}p&\cr
&\dvd{t_{n-2:n-1}}p&&\ \ \ \cdot\ \ \ &&\cr
\dvd{t_{n-1}}p&&\dvd{t_{n-2:n}}p&&&\cr
&\dvd{t_{n-1:n}}p&&&&\cr
\dvd{t_n}p&&&&&\cr
}$$
column by column from left to right, using one of the $n$ pieces of information
$$y(j):= D^{\mu_j}p(t_j),\quad \mu_j:=\#\{i<j: t_i=t_j\}; \quad j=1,\ldots,n,
\label{defdata}
$$
in the leftmost column or else whenever we would otherwise be confronted with
$0/0$.

After construction of this divided difference table, the top diagonal of the
table provides the coefficients
$(\dvd{t_{1:j}}p:j=1,\ldots,n)$ for the Newton form
(with respect to centers $t_1,\ldots,t_{n-1}$)
of the polynomial of degree $<n$ that matches $p$ at $t_{1:n}$, i.e., the
polynomial $p_n$.
More than that, for any sequence $(i_1,\ldots,i_n)$ in which, for each $j$,
$\{i_1,\ldots,i_j\}$ consists of consecutive integers in $[1\fromto n]$,
the above divided difference table provides the coefficients in the Newton form
for the above $r$, but with respect to the centers $(t_{i_j}: j=1{:}n)$.

Now note that the only
information about $p$ entering this calculation is the scalar sequence $y$
described in (\recall{defdata}). Hence we now know the following.

\proclaim Proposition \label{hermiteinterp}.
Let $(t_1,\ldots,t_n)$ have all its multiplicities (if
any) clustered, and let $y\in\FF^n$ be arbitrary. For $j=1,\ldots,n$,
let $c(j)$ be the first entry in the $j$th column in the above divided
difference table as constructed in the described manner from $y$.
\newpar
Then
$$r:=\sum_{j=1}^n\nwt{j-1,t}c(j)$$
is the unique polynomial of degree $<n$ that satisfies the
\dword{Hermite interpolation conditions}
$$ D^{\mu_j}r(t_j)=y(j),\quad \mu_j:=\#\{i<j: t_i=t_j\}; \quad j=1,\ldots,n.
\label{Hermiteconds}
$$

\sect{Evaluation of a Newton form via Horner's method}

\proclaim Horner's method. Let $c(j) := \dvd{t_{1:j}}r$ for
$j=1,\ldots,n>\deg r$, $z\in \FF$, and
$$\eqalign{
\hatc(n)&:=c(n),\cr
\hatc(j)&:=c(j)+(z-t_j)\hatc(j+1),\quad j=n{-}1, n{-}2, \ldots,1.\cr}$$
\newpar
Then
$\hatc(1)=r(z)$.
More than that,
$$
r = \sum_{j=1}^n \nwt{j-1,\hatt\ }\hatc(j),
$$
with
$$
\hatt\;:=\;(z,t_1,t_2,\ldots).
$$

\proof
The first claim follows from the second, or else directly from the fact that
Horner's method is nothing but the evaluation, from the inside out, of
the nested expression
$$
r(z) = c(1) + (z-t_1)(c(2) + \cdots + (z-t_{n-2})(c(n-1) + (z-t_{n-1})c(n))\cdots),
$$
for which reason
Horner's method is also known as \dword{Nested Multiplication}.

As to the second claim, note that $\dvd{z,t_{1:n-1}}r = \dvd{t_{1:n}}r$ since
$\deg r<n$, hence $\hatc(n)=\dvd{z,t_{1:n-1}}r$,
while,
directly from (\recall{compofdvds}),
$$
\dvd{\cdot,t_{1:j-1}} = \dvd{t_{1:j}} + (\cdot-t_j)\dvd{\cdot,t_{1:j}},
\quad j\in\NN,
\label{basicrecur}$$
hence,  by (downward) induction,
$$\hatc(j) = \dvd{z,t_{1:j-1}}r, \quad j=n{-}1,n{-}2,\ldots,1.$$
\endproof


In effect, Horner's Method is another way of filling in a divided difference
table, starting not at the left-most column but with a diagonal, and generating
new entries, not from left to right, but from right to left:

{\ninepoint
$$
\matrix{
\dvd{z}r&&&&\cr
&\dvd{z,t_1}r&&&\cr
\dvd{t_1}r&&\dvd{z,t_{1:2}}r&&\cr
&\dvd{t_{1:2}}r&&\ \ \ \cdot\ \ \ &\cr
\cdot&&\dvd{t_{1:3}}r&&\dvd{z,t_{1:n-2}}r\cr
&\cdot&&\ \ \ \cdot\ \ \ &&\dvd{z,t_{1:n-1}}r=\dvd{t_{1:n}}r\cr
\cdot&&\cdot&&\dvd{t_{1:n-1}}r&\cr
&\cdot&&\ \ \ \cdot\ \ \ &&\dvd{t_{1:n}}r\cr
\cdot&&\ \ \ \cdot\ \ \ &&\cdot&\cr
}$$
}

\noindent
Hence, Horner's method is useful for carrying out a \eword{change of basis},
going from one Newton form to another. Specifically,
$n-1$-fold iteration of this process, with $z=z_{n-1},\ldots,z_1$, is an
efficient way of computing the coefficients $(\dvd{z_{1:j}}r: j=1,\ldots,n)$, of
the Newton form for $r\in\Pi_{<n}$ with respect to the centers $z_{1:n-1}$, from those for
the Newton form with respect to centers $t_{1:n-1}$. Not all the steps need actually be
carried out in case all the $z_j$ are the same, i.e., when switching to the
Taylor form (or local power form).
%
%
%
%
%
%

\sect{Divided differences of functions other than polynomials}

\proclaim Proposition \label{propext}.
On $\Pi$, the divided differences $\dvd{t_{1:j}}$,
$j=1,\ldots,n$, provide a basis for the linear space of linear functionals
spanned by
$$
\dvd{t_j}D^{\mu_j},\quad \mu_j:=\#\{i<j: t_i=t_j\}; \quad j=1,\ldots,n.
\label{listoffnls}
$$

\proof By Proposition \recall{prophermiteconds} and its proof,
$$\cap_{j=1}^n\ker\dvd{t_{1:j}} = \nwt{n,t}\Pi = \cap_{j=1}^n\ker
\dvd{t_j}D^{\mu_j}.$$

Another proof is provided by Horner's method, which, in effect, expresses
$(\dvd{t_{1:j}}: j=1{:}n)$ as linear functions of $(\dvd{t_j}D^{\mu_j}:
j=1{:}n)$, thus showing the first sequence to be contained in the span of the
second. Since the first is linearly independent (as it has $(\nwt{j-1,t}:
j=1{:}n)$ as a dual sequence) while both contain the same number of terms,
it follows that both are bases of the same linear space.
\endproof

This proposition provides a ready extension of $\dvd{t_{1:n}}$ to functions
more general than polynomials, namely to any function for which the derivatives
mentioned in (\recall{listoffnls}) make sense. It is exactly those functions for
which the Hermite conditions (\recall{Hermiteconds}) make sense, hence for which the
Hermite interpolant $r$ of (\recall{hermiteinterp}) is defined. This leads us to G.
Kowalewski's definition.

\proclaim Definition \label{defkow} (\cite{K}). For any smooth enough function $f$
defined, at least, at $t_1,\ldots,t_n$, $\dvd{t_{1:n}}f$ is the leading
coefficient, i.e., the coefficient of $()^{n-1}$, in the power form for the
Hermite interpolant to $f$ at $t_{1:n}$.

In consequence, {\sl $\dvd{t_{1{:}n}}f = \dvd{t_{1{:}n}}p$ for any polynomial
$p$ that matches $f$ at $t_{1{:}n}$}.

\nextex Assume that none of the $t_j$ is zero. Then,
$$\dvd{t_{1:n}}()^{-1} = (-1)^{n-1}/(t_1\cdots t_n).\label{dvdofrecip}$$
This certainly holds for $n=1$
while, for $n>1$,
by (\recall{leibnizformula}), $0=
\dvd{t_{1:n}}(()^{-1}()^1) = \dvd{t_{1:n}}()^{-1}t_n +
\dvd{t_{1:n-1}}()^{-1}$, hence
$\dvd{t_{1:n}}()^{-1} =-\dvd{t_{1:n-1}}()^{-1}/t_n$, and induction finishes the
proof.
This implies the handy formula
$$\dvd{t_{1:n}}(z-\cdot)^{-1} = 1/\nwt{n,t}(z),\quad z\not=t_1,\ldots,t_n.
\label{chakalovone}$$
Therefore, with $\#\xi:=\#\{j: \xi = t_j, 1\le j\le n\}$ the multiplicity with
which $\xi$ occurs in the sequence $t_{1:n}$, and
$$
1/\nwt{n,t}(z) =: \sum_{\xi\in t}\sum_{0\le \mu <\#\xi}{\mu! A_{\xi\mu}\over(z-\xi)^{\mu+1}}
$$
the partial fraction expansion of $1/\nwt{n,t}$, we obtain Chakalov's expansion
$$
\dvd{t_0,\ldots,t_k}f = \sum_{\xi\in t}\sum_{0\le \mu <\#\xi}A_{\xi\mu}D^\mu
f(\xi)
\label{chakalovtwo}
$$
(from \cite{Ch38}) for $f:=1/(z-\cdot)$ for arbitrary $z$ since
$D^\mu 1/(z-\cdot) = \mu!/(z-\cdot)^{\mu+1}$, hence also
for any smooth enough $f$,
by the density of $\{1/(z-\cdot): z\in\FF\}$.

This is an illustration of the peculiar effectiveness of the formula
(\recall{chakalovone}), for
the divided difference of $1/(z-\cdot)$, for deriving and verifying divided
difference identities.
\endexample

\nextex When the $t_j$ are pairwise distinct, (\recall{chakalovtwo}) must
reduce to
$$\dvd{t_{1:n}}f =
\sum_{j=1}^n f(t_j)/D\nwt{n,t}(t_j),
\label{dvdexplicit}$$
since this
is readily seen to be the leading coefficient of the
polynomial of degree $<n$ that matches a given $f$ at the $n$
pairwise distinct sites $t_1,\ldots,t_n$ when we write that polynomial
in Lagrange form,
$$\sum_{j=1}^n f(t_j)\prod_{i\in 1{:}n\bs j}{\cdot-t_i\over t_j-t_i}.$$

It follows (see the proof of \cite{ET:Lemma I})
that, {\sl for $-1\le t_1<\cdots<t_n\le1$,
$$
\norm{\dvd{t_{1:n}}:C([-1\fromto1])\to\FF}
= \sum_{j=1}^n 1/|D\nwt{n,t}(t_j)| \ge 2^{n-2},
\label{lowerbound}
$$
with equality iff $\nwt{n,t} = (()^2-1)U_{n-2}$}, where $U_{n-2}$ is
the second-kind Chebyshev polynomial.

Indeed, for any such
$$\gt:=(t_1,\ldots,t_n),$$
the restriction
$\gl$ of $\dvd{\gt}$ to $\Piln$ is the
unique linear functional on $\Piln$ that vanishes on $\Pi_{<n-1}$ and takes the
value $1$ at $()^{n-1}$, hence takes its norm on the error of the best
(uniform) approximation to $()^{n-1}$ from $\Pi_{<n-1}$, i.e., on the Chebyshev
polynomial of degree $n-1$. Each such $\dvd{\gt}$ is an extension of this $\gl$,
hence has norm $\ge\norm{\gl} = 1/\dist(()^{n-1},\Pi_{<n-1}) = 2^{n-2}$, with
equality iff $\dvd{\gt}$ takes on its norm on that Chebyshev polynomial, i.e.,
iff $\gt$ is the sequence of extreme sites of that Chebyshev polynomial.
\endexample

\sect{The divided difference as approximate normalized derivative}
Assume that $f$ is differentiable on an interval that contains the
nondecreasing finite sequence
$$\gt = (\gt_0\le\cdots\le\gt_k),$$
and assume further that $\dvd{\gt}f$ is defined, hence so is the Hermite
interpolant
$$P_\gt f$$
of $f$ at $\gt$.

Then $f-P_\gt f$ vanishes at $\gt_{0{:}k}$, therefore $D(f-P_\gt f)$
vanishes at some $\gs = (\gs_0,\ldots,\gs_{k-1})$ that \dword{interlaces} $\gt$,
meaning that
$$\gt_i\le\gs_i\le\gt_{i+1},\quad {\rm all\ } i.$$
This is evident when $\gt_i=\gt_j$ for some $i<j$, and is Rolle's Theorem when
$\gt_i<\gt_{i+1}$.
Consequently, $DP_\gt f$ is a polynomial of degree $<k$ that matches $Df$ at
$\gs_0,\ldots,\gs_{k-1}$, hence must be its Hermite interpolant at $\gs$.
This proves the following.

\proclaim Proposition \label{propderiv} (\cite{Ho}).
If $f$ is differentiable on an interval that contains the nondecreasing
$(k+1)$-sequence $\gt$ and smooth enough at $\gt$ so that its Hermite
interpolant, $P_\gt f$, at $\gt$ exists, then there is a nondecreasing
$k$-sequence $\gs$ interlacing $\gt$ and so that
$$P_\gs (Df) = D P_\gt f.$$
In particular, then
$$k\dvd{\gt}f = \dvd{\gs}Df.$$

From this, induction provides the
\proclaim Corollary (\cite{Sc}). Under the same assumptions, but with $f$
$k$ times differentiable on that interval, there exists $\xi$ in that
interval for which
$$k!\dvd{\gt_0,\ldots,\gt_k} = D^kf(\xi).
\label{intermedval}$$

The special case $k=1$, i.e.,
$$\dvd{a,b}f = Df(\xi), \quad \hbox{for\ some\ } \xi\in (a\fromto b),$$
is so obvious a
consequence or restatement of L'H\^opital's Rule, it must have been around at
least that long.

Chakalov \cite{Ch34} has made a detailed study of the possible values that $\xi$
might take in (\recall{intermedval}) as $f$ varies over a given class of functions.

\cite{KA:p.~91} reports that already Taylor, in \cite{T}, derived his
eponymous expansion (\recall{taylorexp}) as the limit of Newton's formula, albeit for
equally spaced sites only.

\sect{Representations}

\fsubsect{Determinant ratio}
Let
$$\gt := (\gt_0,\ldots,\gt_k).$$
Kowalewski's definition of $\dvd{\gt}f$ as the leading coefficient, in the power
form, of the Hermite interpolant to $f$ at $\gt$ gives,
for the case of simple sites and via Cramer's Rule, the formula
$$
\dvd{\gt}f\;=\;\det Q_\gt[()^0,\ldots,()^{k-1},f]/\det Q_\gt[()^0,\ldots,()^k]
\label{eqdet}
$$
in which
$$Q_\gt [g_0,\ldots,g_k] := (g_j(\gt_i): i,j=0,\ldots,k).$$
In some papers and books, the identity (\recall{eqdet}) serves as the definition
of $\dvd{\gt}f$ despite the fact that it needs awkward modification in
the case of repeated sites.

\subsect{Peano kernel (B-spline)}
Assume that $\gt:= (\gt_0,\ldots,\gt_k)$ lies in the interval $[a\fromto b]$ and that $f$ has $k$
derivatives on that interval. Then, on that interval, we have
Taylor's identity
$$
f(x) = \sum_{j<k}(x-a)^jD^jf(a)/j! + \int_a^b(x-y)_+^{k-1} D^kf(y)\dd y/(k-1)!.
\label{eqtaylor}$$
If now $\gt_0<\gt_k$, then, from Proposition \recall{propext}, $\dvd{\gt}$ is a
weighted sum of values of derivatives of order $<k$, hence commutes with the
integral in Taylor's formula (\recall{eqtaylor}) while, in any case, it annihilates any
polynomial of degree $<k$. Therefore
$$
\dvd{\gt}f = \int_a^b M(\cdot|\gt) D^kf/k!,
\label{eqPeano}
$$
with
$$M(x|\gt):= k\dvd{\gt}(\cdot-x)_+^{k-1}
\label{eqBspline}
$$
the Curry-Schoenberg \dword{B-spline} (see \cite{CS}) with knots $\gt$ and
normalized to have integral 1.
While Schoenberg and Curry named and studied the B-spline only in the 1940's,
it appears in this role as the Peano kernel for the divided difference already
earlier, e.g., in \cite{P33} and \cite{Ch34} (see \cite{BP}) or \cite{Fa}.

\subsect{Contour integral}
An entirely different approach to divided differences and Hermite interpolation
begins with Frobenius' paper \cite{Fr}, so different that it had no influence
on the literature on interpolation (except for a footnote-like mention in
\cite{Ch38}). To be sure, Frobenius himself seems to have
thought of it more as an exercise in expansions, never mentioning the word
`interpolation'. Nevertheless, Frobenius describes in full detail the salient
facts of polynomial interpolation in the complex case, with the aid of the
Cauchy integral.

In \cite{Fr}, Frobenius investigates \dword{Newton series}, i.e.,
\eword{infinite} expansions
$$\sum_{j=1}^\infty c_j \nwt{j-1,t}$$
in  the Newton polynomials
$\nwt{j,t}$ defined in (\recall{defnwt}).
He begins with the identity
$$(y-x)\sum_{j=1}^n{\nwt{j-1,t}(x)\over\nwt{j,t}(y)}
\;=\;1 - {\nwt{n,t}(x)\over\nwt{n,t}(y)},\label{frobeniusid}$$
a ready consequence of the observations
$$\eqalign{x\nwt{j-1,t}(x)\; & \;=\; \nwt{j,t}(x)\;\; +\;\;
t_{j}\nwt{j-1,t}(x),\cr
\noalign{\vskip2pt}
{y\over\nwt{j,t}(y)}\; & \;=\;\; {1\over\nwt{j-1,t}(y)} \;\;+\;\;{t_{j}\over\nwt{j-1,t}(y)}\cr}$$
since these imply that
$$
\eqalign{
y\sum_{j=1}^n{\nwt{j-1,t}(x)\over\nwt{j,t}(y)}
\;=\; &\sum_j{\nwt{j-1,t}(x)\over\nwt{j-1,t}(y)}
+ \sum_j{t_{j}\nwt{j-1,t}(x)\over\nwt{j,t}(y)},\cr
\noalign{\vskip2pt}
x\sum_{j=1}^n{\nwt{j-1,t}(x)\over\nwt{j,t}(y)}
\;=\; &\sum_j{\nwt{j,t}(x)\over\nwt{j,t}(y)}
+ \sum_j{t_{j}\nwt{j-1,t}(x)\over\nwt{j,t}(y)}.\cr}
$$
Then (in \S 4), he uses (\recall{frobeniusid}), in the form
$$
\sum_{j=1}^n{\nwt{j-1,t}(z)\over\nwt{j,t}(\gz)}
 + {\nwt{n,t}(z)\over(\gz-z)\nwt{n,t}(\gz)}
\;=\;1/(\gz-z),
$$
in Cauchy's formula
$$f(z) = \ootpii\oint {f(\gz)\dd \gz\over \gz-z}
$$
to conclude that
$$f(z) = \sum_{j=1}^n\nwt{j-1,t}c_j\;+\;
\nwt{n,t}\ootpii\oint{f(\gz)\dd\gz\over(\gz-z)\nwt{n,t}(\gz)},
\label{frobenius}$$
with
$$c_j\;:=\;\ootpii\oint{f(\gz)\dd\gz\over\nwt{j,t}(\gz)},\quad j=1,\ldots,n.$$
For this, he assumes that $z$ is in some disk of radius $\rho$, in which $f$ is
entire, and $\gz$ runs on the boundary of a disk of radius $\rho'<\rho$ that
contains $z$, with none of the relevant $t_j$ in the annulus formed by
the two disks.

Directly from the definition of the divided difference, we therefore conclude
that, under these assumptions on $f$ and $t$,
$$\dvd{t_{1{:}j}}f
\;=\;\ootpii\oint{f(\gz)\dd\gz\over\nwt{j,t}(\gz)},\quad j=0,1,2, \ldots.
\label{dvdaslineintegral}$$
Strikingly, Frobenius never mentions that (\recall{frobenius}) provides a
general polynomial interpolant and its error. Could he have been unaware of it?
To be sure, he could not have called it `Hermite interpolation' since Hermite's
paper \cite{H78} appeared well after Frobenius'. There is no indication that
Hermite was aware of Frobenius' paper.

\subsect{Genocchi-Hermite}
Starting with (\recall{dvdisdvd}) and the observation that
$$\dvd{x,y}f = \int_0^1 Df((1-s)x + sy)\dd s,$$
induction on $n$ gives the
(univariate) Genocchi-Hermite formula
$$
\dvd{\gt_0, \ldots, \gt_n}f = \int_{[\gt_0, \ldots, \gt_n]} D^nf,
\label{genocchi}
$$
with
$$
\eqalign{\int_{[\gt_0, \ldots, \gt_n]} f\;:=\;&\cr
 &\kern-2truecm\int_0^1\int_0^{s_1}\cdots\int_0^{s_{n-1}}
f((1-s_1)\gt_0 +  \cdots + (s_{n-1}-s_n)\gt_{n-1} + s_n\gt_n)
\dd s_n\cdots \dd s_1.\cr}
$$

\cite{No:p.16} mistakenly attributes (\recall{genocchi})
to \cite{H59}, possibly because that paper carries the suggestive title
``Sur l'interpolation''.

At the end of the paper \cite{H78}, on polynomial interpolation
to data at the $n$ pairwise distinct sites $t_1,\ldots,t_n$ in the complex
plane, Hermite does give a formula involving the righthand-side of the above,
namely the formula
$$
f(x) - Pf(x) = (x-t_1)\cdots(x-t_n)\int_{[t_n,\ldots,t_1,x]} D^nf
$$
for the error in the Lagrange interpolant $Pf$ to $f$ at $t_{1:n}$.
Thus, it requires the observation that
$$
f(x) - Pf(x) = (x-t_1)\cdots(x-t_n)\dvd{t_n,\ldots,t_1,x}f
$$
to deduce the Genocchi-Hermite formula from \cite{H78}.
(He also gives the rather more complicated formula
$$
\eqalign{f(x)-Pf(x) =&\cr
&\kern-2truecm (x-a_1)^\ga\cdots(x-a_n)^\gl
\int_{[a_n,\ldots,a_1,x]}
\braket{s_n-s_{n-1}}^{\ga-1}\cdots\braket{1-s_1}^{\gl-1}D^{\ga+\cdots+\gl}f\cr}
$$
 for the error in case of repeated interpolation. Here, $\braket{z}^j:=z^j/j!$.)

In contrast, \cite{G69} is explicitly concerned with a
\eword{representation formula for the divided difference}.
However, the `divided difference' he represents is the following:
$$
\dvd{x,x+h_1}\dvd{\cdot,\cdot+h_2}\cdots\dvd{\cdot,\cdot+h_n}
=(\Delta_{h_1}/h_1)\cdots(\Delta_{h_n}/h_n)
$$
and for it he gets the representation
$$
\int_0^1\cdots\int_0^1
D^n f(x + h_1t_1 + \cdots + h_nt_n)
\dd t_1\cdots \dd t_n.
$$

\cite{No:p.16} cites \cite{G78a}, \cite{G78b} as places where formulations
equivalent to the Genocchi-Hermite formula can be found. So far, I've been
only able to find \cite{G78b}. It is a letter to Hermite, in which
Genocchi brings, among other things, the above representation formula to
Hermite's attention, refers to a paper of his in [Archives de Grunert, t. XLIX,
3e cahier] as containing a corresponding error formula for Newton
interpolation. He states that he, in continuing work, had obtained such a
representation also for Amp\`ere's fonctions interpolatoires (aka divided
differences), and finishes with the formula
$$
\eqalign{
\int_0^1\cdot\cdot\int_0^1
s_1^{n-1} s_2^{n-2} \cdots s_{n-1}&\cr
&\kern-2.5truecmD^nf (x_0  + s_1(x_1-x_0) +  \cdots
+ s_1 s_2 \cdots s_n (x_n-x_{n-1}))
\dd s_1\cdots \dd s_n\cr}
$$
for $\dvd{x_0,\ldots,x_n}f$, and says that it is equivalent to the formula
$$
\dvd{x_0,\ldots,x_n}f =
\int\cdots\int
D^nf (s_0x_0  + s_1x_1 + \cdots s_nx_n)
 \dd s_1\cdots \dd s_n
$$
in which the conditions $s_0+\cdots+s_n=1$, $s_i\ge0$, all $i$, are imposed.

\cite{St27:p.17f} proves the Genocchi-Hermite formula but calls it
Jensen's formula, because of \cite{J}.

\sect{Divided difference expansions of the divided difference}
By applying $\dvd{s_{1:m}}$ to both sides of the identity
$$
\dvd{\cdot} = \sum_{j=1}^n \nwt{j-1,t}\dvd{t_{1:j}}
+ \nwt{n,t}\dvd{t_{1:n},\cdot}$$
obtained from (\recall{withremainder}),
one obtains the expansion
$$
\dvd{s_{1:m}} =
\sum_{j=m}^n \dvd{s_{1:m}}\nwt{j-1,t}\dvd{t_{1:j}} +
E(s,t),$$
where, by the Leibniz formula (\recall{leibnizformula}),
$$
\eqalign{
E(s,t) :=\;&
\dvd{s_{1:m}}(\nwt{n,t} \dvd{t_{1:n},\cdot})\cr
=\;&\sum_{i=1}^m \dvd{s_{i:m}}\nwt{n,t} \dvd{s_{1:m},t_{1:n}}.
}
$$

But, following \cite{Fl} and with $p:=n-m$, one gets the better formula
$$
E(s,t) :=
\sum_{i=1}^m(s_i-t_{i+p})
(\dvd{s_{1:i}}\nwt{i+p,t}) \dvd{t_{1:i+p},s_{i:m}}
$$
in which all the divided differences on the right side are of the same order, $n$.
The proof (see \cite{B03}), by induction on $n$, uses the easy consequence of
(\recall{leibnizformula}) that
$$
(s_i-y)\dvd{s_{i:m}}f =
\dvd{s_{i:m}}((\cdot-y)f)-\dvd{s_{i+1:m}}f.$$
The induction is anchored at $n=m$ for which the formula
$$
\dvd{s_{1:m}} - \dvd{t_{1:m}}
=\sum_{i=1}^m (s_i-t_i)\dvd{s_{1:i},t_{i:m}}
$$
can already be found in \cite{Ho}.

\sect{Notation and nomenclature}
It is quite common in earlier literature to use the notation
$$[y_1,\ldots,y_j]$$
for the divided difference of order $j-1$ of data $((t_i,y_i): i=1{:}j)$.
This reflects the fact that divided differences were thought of as convenient
expressions in terms of the given data rather than as linear functionals on
some vector space of functions.

The presently most common notation for $\dvd{t_{1:j}}p = \dvd{t_1,\ldots,t_j}p$
is
$$p[t_1,\ldots,t_j]$$
(or, perhaps, $p(t_1,\ldots,t_j)$)
which enlarges upon the fact that $\dvd{z}p = p(z)$,  but this becomes awkward
when the divided difference is to be treated as a linear functional.
In that regard, the notation
$$[t_1,\ldots,t_j]p$$
is better, but suffers from the fact that the resulting notation
$$[t_1,\ldots,t_j]$$
for the linear functional itself conflicts with standard notations, such as the
matrix (or, more generally, the column map) with columns $t_1,\ldots, t_j$, or,
in the special case $j=2$, i.e.,
$$[t_1,t_2],$$
the closed interval with endpoints $t_1$ and $t_2$ or else the scalar product
of the vectors $t_1$ and $t_2$. The notation
$$[t_1,\ldots,t_j;p]$$
does not suffer from this defect, as it leads to the notation
$[t_1,\ldots,t_j;\cdot]$ for the linear functional itself, though it requires
the reader not to mistakenly read that semicolon as yet another comma.

The notation, $\divdif$, used in this essay was proposed by W. Kahan
some time ago (see, e.g., \cite{Ka}), and does not suffer from any of the
defects mentioned and has the advantage of being literal (given that $\gD$
is standard notation for a difference). Here is a \TeX\ macro for it:

\medskip
\def\bs{\char'134}
\def\lb{\char'173}
\def\rb{\char'175}
\centerline{\vbox{\tt\noindent
\bs def\bs divdif\lb\bs mathord{\bs kern.43em\bs vrule width.6pt height5.6pt
\newline\phantom{\bs def\bs divdif\lb\bs}
depth-.28pt \bs kern-.43em\bs Delta\rb}
}}
\medskip

Although divided differences are rightly associated with Newton (because of
\cite{N87:Book iii, Lemma v, Case ii}, \cite{N11}; for a detailed discussion,
including facsimiles of the originals and their translations, see
\cite{Fra18}, \cite{Fra19}, \cite{Fra27}),
the term `divided difference' was,
according to \cite{WR:p.20},
first used in \cite{M:p.550}, even though, by then,
Amp\`ere \cite{A} had called it \dword{fonction
interpolaire}, and this is the term used in the French literature of the 1800s.
To be sure, in \cite{N11}, Newton repeatedly uses the term ``differentia
sic divisa'.




\def\AM{Appl.\ Math.}
\def\CRASP{C. R. Acad.\ Sci.\ Paris}
\def\JAM{J. Analyse Math.}
\def\JAT{J. Approx.\ Theory}
\def\JRAM{J. reine angew.\ Math.}
\def\SJNA{SIAM J. Numer.\ Anal.}
\def\ZAMM{Z. Angew.\ Math.\ Mech.}
\References

\rhl{A:Amp\`ere 1826}
\refJ Amp\`ere, Andr\'e-Marie;
Essai sur un noveau mode d'exposition des princi\-pes du calcul diff\'erentiel,
  du calcul aux diff\'erences et de l'interpolation des suites,
    consid\'er\'ees comme d\'erivant d'un source commune;
Annals de mathematicques pures et appliqu\'ees de Gergonne (Paris, Bachelier,
  in $4^o$); 16; 1826; 329--349;

\rhl{B72:de Boor 1972}
\refJ de Boor,  C.;
On calculating with $B$-splines;
\JAT; 6; 1972; 50--62;

\rhl{B03:de Boor 2003a}
\refJ de Boor, C.;
A divided difference expansion of a divided difference;
\JAT; 122; 2003a; 10--12;

\rhl{B03b:de Boor 2003b}
\refJ de Boor, C.;
A Leibniz formula for multivariate divided differences;
\SJNA; 41(3); 2003b; 856--868;

\rhl{BP:de Boor et al. 2003}
\refJ de Boor, C., Pinkus, A.;
The B-spline recurrence relations of Chakalov and of Popoviciu;
\JAT; 124; 2003; 115--123;

\rhl{BR:Bulirsch et al. 1968}
\refQ Bulirsch,  R., Rutishauser, H.;
Interpolation und gen\"aherte Quadratur;
(Mathematische Hilfsmittel des Ingenieurs, Teil III), S. Sauer \& I. Szab\'o
(eds.), Grundlehren vol.~141, Springer-Verlag (Heidelberg); 1968; 232--319;

\rhl{Ca:Cauchy 1840}
\refJ Cauchy, Augustin;
Sur les fonctions interpolaires;
\CRASP; 11; 1840; 775--789;

\rhl{Ch34:Tchakaloff 1934}
\refJ Tchakaloff,  L.;
Sur la structure des ensembles lin\'eaires d\'efinis par une certaine
   propri\'et\'e minimale;
Acta Math; 63; 1934; 77--97;

\rhl{Ch38:Chakalov 1938}
\refJ Chakalov,  L.;
On a certain presentation of the Newton divided differences in interpolation
   theory and its applications (in Bulgarian);
Annuaire Univ.\ Sofia, Fiz.\ Mat.\ Fakultet; 34; 1938; 353--394;

\rhl{CS:Curry \& Schoenberg 1966}
\refJ Curry,  H. B., Schoenberg, I. J.;
On P\'olya frequency functions IV: the fundamental spline
     functions and their limits;
\JAM; 17; 1966;  71--107;

\rhl{ET:Erd\H os et al. 1940}
\refJ Erd\H os, P., Tur\'an, P.;
On interpolation, III. Interpolatory theory of polynomials;
\AM (2); 41; 1940; 510--553;

\rhl{Fa:Favard 1940}
\refJ Favard,  J.;
Sur l'interpolation;
J.  Math.\  Pures Appl.; 19; 1940; 281--306;

\rhl{Fl:Floater 2003}
\refJ Floater, M.;
Error formulas for divided difference expansions and numerical
   differentiation;
\JAT; 122; 2003; 1--9;

\rhl{Fra18:Fraser 1918}
\refJ Fraser, Duncan C.;
Newton's interpolation formulas;
J. Instit.\ Actuaries; LI; 1918; 77--106;

\refJ Fraser, Duncan C.;
\rhl{Fra19:Fraser 1919}
Newton's interpolation formulas. Further Notes;
J. Instit.\ Actuaries; LI; 1919; 211--232;

\rhl{Fra27:Fraser 1927}
\refJ Fraser, Duncan C.;
Newton's interpolation formulas. An unpublished manuscript of Sir Isaac
   Newton;
J. Instit.\ Actuaries; LVIII; 1927; 53--95;

\rhl{Fr:Frobenius 1871}
\refJ Frobenius, G.;
\"Uber die Entwicklung analytischer Functionen in Reihen, die nach
gegebenen Functionen fortschreiten;
\JRAM; 73; 1871; 1--30;

\rhl{G69:Genocchi 1869}
\refJ Genocchi, A.;
Relation entre la diff\'erence et la d\'eriv\'ee d'un m\^eme ordre quelconque;
Archiv Math.\ Phys.\ (I); 49; 1869; 342--345;

\rhl{G78a:Genocchi 1878a}
\refJ Genocchi, A.;
Intorno alle funzioni interpolari;
Atti della Reale Accademia delle Scienze di Torino; 13; 1878a; 716--729;

\rhl{G78b:Genocchi 1878b}
\refJ Genocchi, A.;
Sur la formule sommatoire de Maclaurin et les fonctions interpolaires;
\CRASP; 86; 1878b; 466--469;

\rhl{H59:Hermite 1859}
\refJ Hermite, Ch.;
Sur l'interpolation;
\CRASP; 48; 1859; 62--67;

\rhl{H78:Hermite 1878}
\refJ Hermite, Ch.;
Formule d'interpolation de Lagrange;
\JRAM; 84; 1878; 70--79;

\rhl{Ho:Hopf 1926}
\refD Hopf, Eberhard;
\"Uber die Zusammenh\"ange zwischen gewissen h\"oheren
   Diffe\-renzen\-quo\-tienten reeller Funktionen einer reellen Variablen und
  deren Dif\-fe\-renzier\-bar\-keits\-eigen\-schaften;
Universit\"at Berlin (30 pp); 1926;

\rhl{J:Jensen 1894}
\refJ Jensen, J. L. W. V.;
Sure une expression simple du reste dans la formule d'interpolation de
Newton;
Bull.\ Acad.\ Roy.\ Danemark; xx; 1894; 246--xxx;

\rhl{Ka:Kahan 1974}
\refR Kahan, W.;
Divided differences of algebraic functions;
Class notes for Math 228A, UC Berkeley, Fall; 1974;

\rhl{KA:A. Kowalewski 1917}
\refB Kowalewski, Arnold;
Newton, Cotes, Gauss, Jacobi: Vier grund\-legende Ab\-hand\-lungen \"uber
   Interpolation und gen\"aherte Quadratur;
Teubner (Leipzig); 1917;

\rhl{K:G. Kowalewski 1932}
\refB Kowalewski, G.;
Interpolation und gen\"aherte Quadratur;
B. G. Teubner (Berlin); 1932;

\rhl{M:de Morgan 1842}
\refB de Morgan, A.;
The Differential and Integral Calculus;
Baldwin and Cradock (London, UK); 1842;

\rhl{N76:Newton 1676}
\refR Newton, I.;
Epistola posterior (Second letter to Leibniz via Oldenburg);
24 oct; 1676;

\rhl{N87:Newton 1687}
\refB Newton, I.;
Philosophi\ae\ naturalis principia mathematica;
Joseph Strea\-ter (London); 1687;
(see \newline{\tt ftp://ftp.math.technion.ac.il/hat/fpapers/newton1.pdf}
 and\newline
{\tt ftp://ftp.math.technion.ac.il/hat/fpapers/newton2.pdf})

\rhl{N11:Newton 1711}
\refB Newton, I.;
Methodus differentialis;
Jones (London); 1711;

\rhl{No:N\"orlund 1924}
\refB N\"orlund,  N. E.;
Vorlesungen \"uber Differenzenrechnung;
Grundlehren Vol.\ XIII, Springer (Berlin); 1924;

\rhl{O:Opitz 1964}
\refJ Opitz, G.;
Steigungsmatrizen;
\ZAMM; 44; 1964; T52--T54;

\rhl{P33:Popoviciu 1933}
\refD Popoviciu, Tib{\`e}re;
Sur quelques propri\'et\'es des fonctions d'une ou de deux variables
r\'eelles;
presented to the Facult\'e des Sciences de Paris, published by Institutul de
Arte Grafice ``Ardealul'' (Cluj, Romania); 1933;

\rhl{P40:Popoviciu 1940}
\refJ Popoviciu, Tiberiu;
Introduction \`a la th\'eorie des diff\'erences divis\'ees;
Bull.\ Mathem., Societea Romana de Stiinte, Bukharest; 42; 1940; 65--78;

\rhl{Sc:Schwarz 1881-2}
\refJ Schwarz, H. A.;
D\'emonstration \'el\'ementaire d'une propri\'et\'e fondamentale des fonctions
   interpolaires;
Atti Accad.\ Sci.\ Torino; 17; 1881-2; 740--742;

\rhl{St27:Steffensen 1927}
\refB Steffensen, J. F.;
Interpolation;
Williams and Wilkins (Baltimore); 1927;

\rhl{St39:Steffensen 1939}
\refJ Steffensen, J. F.;
Note on divided differences;
Danske Vid.\ Selsk.\ Math.-Fys.\ Medd.; 17(3); 1939; 1--12;

%
\rhl{T:Taylor 1715}
\refB Taylor, Brook;
Methodus incrementorum directa et inversa;
London (England); 1715;

\rhl{WR:Whittaker et al. 1937}
\refB Whittaker, E. T., Robinson, G.;
The Calculus of Observations, A Treatise on Numerical Mathematics, 2nd ed.;
Blackie (Glasgow, UK); 1937;



{

\bigskip\obeylines
C. de Boor
P.O. Box 1076
Eastsound WA 98245 USA
{\tt deboor@cs.wisc.edu}
{\tt http://www.cs.wisc.edu/}$\sim${\tt deboor/}

}


\bye